\documentclass{amsart}

\usepackage{latexsym}
\usepackage{amssymb}
\usepackage{amsmath,mathabx,mathtools}
\usepackage{graphicx}

\newtheorem{theorem}{Theorem}[section]
\newtheorem{lemma}[theorem]{Lemma}

\theoremstyle{definition}
\newtheorem{definition}[theorem]{Definition}

\theoremstyle{remark}

\numberwithin{equation}{section}

\newcommand{\R}{\mathbb{R}}
\newcommand{\Z}{\mathbb{Z}}

\newcommand{\lm}{\lambda}
\newcommand{\Lm}{\Lambda}
\newcommand{\Ga}{\Gamma}

\newcommand{\f}{\varphi}
\newcommand{\SL}{\mathrm{SL}}
\newcommand{\PSL}{\mathrm{PSL}}

\newcommand{\B}{\mathbb B}

\newcommand{\Pp}{\mathbb P}

\newcommand{\Pb}{\mathrm{Prob}}
\newcommand{\Ss}{\mathcal S}

\newcommand{\K}{\mathbb K}

\newcommand{\la}{\langle}
\newcommand{\ra}{\rangle}

\begin{document}

\title[Linear actions and solid factors]{Solid factors associated to the linear action of $\SL_2(\Z)$ on $\Z^n$}

\author{Paul Jolissaint}
\address{Universit\'e de Neuch\^atel,
       Institut de Math\'ematiques,       
       E.-Argand 11,
       2000 Neuch\^atel, Switzerland}
       
\email{paul.jolissaint@ik.me}


\date{\today}


\begin{abstract}
Let $n\geq 2$ be an integer. Set $\Ga_n=\SL_2(\Z)$ (respectively $\PSL_2(\Z)$) for $n$ even (respectively for $n$ odd). Then $\Ga_n$ acts on $\Z^n$ through the restriction of the linear action of $\SL_2(\R)$ on $\R^n$. We prove that the von Neumann algebra $L(\Z^n\rtimes \Ga_n)$ is a solid factor. This generalizes the case $n=2$ proved by N. Ozawa in \cite{Oz1}.
\end{abstract} 

\maketitle

\section{Introduction}

Let $n\geq 2$ be an integer and let $P_{n-1}(\Z)$ be the ring of homogeneous polynomials of degree $n-1$ in two variables $X,Y$ with integer coefficients: we have a standard action $\rho_{n-1}$ of $\SL_2(\Z)$ on $P_{n-1}(\Z)$ defined by 
\[
\Big(\rho_{n-1}\begin{pmatrix}
a & b \\c & d\end{pmatrix}\Big)(P)(X,Y)=P(aX+cY,bX+dY)
\]
for $P\in P_{n-1}(\Z)$ and 
$\begin{pmatrix}
a & b \\c & d
\end{pmatrix}
\in \SL_2(\Z)$. Identifying the abelian group $P_{n-1}(\Z)$ with $\Z^n$ by means of the canonical basis $X^{n-1},X^{n-2}Y,\ldots,XY^{n-2},Y^{n-1}$, we get an $\SL_2(\Z)$-action on $P_{n-1}(\Z)=\Z^n$. If $n$ is odd, the representation $\rho_{n-1}$ factors through $\PSL_2(\Z)$. Setting $\Ga_n=\SL_2(\Z)$ for $n$ even and $\Ga_n=\PSL_2(\Z)$ for $n$ odd, the associated group von Neumann algebra $L(\Z^n\rtimes \Ga_n)$ is a type II$_1$ factor, as it was proved in \cite{JV}.

\medskip

 The present note consists in a generalization of the main result of the short note \cite{Oz1} where the author proves that the group von Neumann algebra $L(\Z^2\rtimes \SL_2(\Z))$ is \textit{solid}, i.e. the relative commutant of every diffuse von Neumann subalgebra of $L(\Z^2\rtimes \SL_2(\Z))$ is injective. 

Solid von Neumann algebras were introduced and studied first in the seminal article \cite{Oz2} as a sufficient condition which ensures that a solid type II$_1$ factor is necessarily non-$\Ga$ and prime unless it is injective. Thus, using ideas and techniques of \cite{Oz2,Oz3,Oz1}, our main result is the following theorem.

\begin{theorem}
Let $n\geq 2$ be an integer. Then the type $\mathrm{II}_1$ factor $L(\Z^n\rtimes \Ga_n)$ is solid.
\end{theorem}

As will be explained in Section 2, the above theorem is a consequence of the following slightly more general result. We refer to Section 2 and to \cite{BO} for precise definitions.

\begin{theorem}
Let $\Ga$ and $\Lm$ be countable groups. Assume that:
\begin{enumerate}
\item [(1)] $\Ga$ is a bi-exact group that acts on $\Lm$;
\item [(2)] $\Lm$ is an amenable group;
\item [(3)] $\Ga$ acts amenably on a metric space $(X,d)$ and there is a $\Ga$-equivariant map $\f:\Lm\setminus\{e\}\rightarrow X$ such that
\[
\lim_{t\to\infty}d(\f(sts'),\f(t))=0
\]
for all $s,s'\in \Lm$.
\end{enumerate}
Then the von Neumann algebra $L(\Lm\rtimes \Ga)$ is solid.
\end{theorem}

\medskip
\textit{Acknowledgements.} I am very grateful to A. Valette for his help in Lemma 2.4.

\section{Proofs of Theorems 1.1 and 1.2}

Let $G$ be a discrete, countable group. If $G$ acts on some set $S$ and if $f$ is any function defined on $S$ then we denote by $g.f$ the function $s\mapsto f(g^{-1}.s), g\in G$.
\medskip

Recall from \cite[Theorem 5.1.7]{BO} that $G$ is \textit{exact} if it admits an amenable action on some compact topological space, i.e. if there is a compact space $X$ on which $G$ acts by homeomorphisms so that there exists a sequence of continuous maps $m_i:X\rightarrow \Pb(G)$ such that, for every $g\in G$
\[
\lim_{i\to\infty}\Big(\sup_{x\in X}\Vert g.m_i^x-m_i^{g.x}\Vert_1\Big)=0.
\]

\medskip\noindent
\textbf{Notations} (1) Let $A,B$ be $C^*$-algebras. We denote by $A\odot B$ the algebraic tensor product of $A$ and $B$, and by $A\otimes B$ their minimal tensor product.\\
(2) If $H$ is a Hilbert space, $\K(H)$ denotes the ideal of all compact operators in the algebra of all linear, bounded operators $\B(H)$.\\
(3) Let $G$ be a group. We denote by $\lm$ (resp. $\rho$) its left (resp. right) regular representation on $\ell^2(G)$, and by $C^*_\lm(G)$ and $C^*_\rho(G)$ the $C^*$-algebras generated by $\lm(G)$ and $\rho(G)$ respectively. We write $\lm_G$ and $\rho_G$ when the group must be precised.
Finally, we denote by $L(G)$ the associated von Neumann algebra of $G$, which is the $\sigma$-weak closure of the reduced $C^*$-algebra $C^*_\lm(G)$.\\

\begin{definition}
\cite[Definition 15.1.2]{BO}
The countable group $G$ is \textit{bi-exact} (or has \textit{property} $\Ss$) if it is exact and there is a map $\mu:G\rightarrow \Pb(G)$ such that, for all $g,h\in G$, 
\[
\lim_{x\to\infty}\Vert \mu(gxh)-g.\mu(x)\Vert_1=0.
\]
\end{definition}

\medskip
We state now the following criterion of bi-exactness for exact groups. It is taken from \cite[Lemma 15.1.4]{BO}, and the last assertion is \cite[Theorem 6]{Oz2}.

\begin{theorem}
Let $G$ be a countable exact group. Then the following conditions are equivalent:
\begin{enumerate}
\item [(1)] $G$ is bi-exact;
\item [(2)]  there exists a unital completely positive (ucp) map 
\[
\Phi:C^*_\lm(G)\otimes C^*_\rho(G)\rightarrow \B(\ell^2(G))
\] 
such that $\Phi(a\otimes b)-ab\in \K(\ell^2(G))$ for all $a\in C^*_\lm(G),b\in C^*_\rho(G)$.
\end{enumerate}
If it is the case, then $L(G)$ is a solid von Neumann algebra: for every diffuse von Neumann subalgebra $N\subset L(G)$, its relative commutant $N'\cap L(G)$ is injective.
\end{theorem}

Condition (2) above was first established for the free group on two generators by Akemann and Ostrand in \cite{AO}, later by Skandalis \cite{Sk} for discrete subgroups of connected simple Lie groups of rank one, and by Higson and Guentner for hyperbolic groups in \cite{HG}. In particular, $\SL_2(\Z)$ and $\PSL_2(\Z)$ are bi-exact.

As we will prove, the semi-direct products $\Z^n\rtimes \Ga_n$ defined in Section 1 are bi-exact because they possess an extra property (introduced in \cite{Oz1} for the case $n=2$) that implies their bi-exactness. 

\medskip
From now on, let $\Ga$ be a countable, bi-exact group which acts on a countable, amenable group $\Lm$ via an action denoted by $\alpha$. We assume moreover that the pair $(\Ga,\Lm)$ satisfies the following additional condition: there is a compact metric space $(X,d)$ on which $\Ga$ acts amenably and there is a $\Ga$-equivariant map $\f:\Lm\setminus \{e\}\rightarrow X$ such that, for all $s,s'\in \Lm$,
\[
\lim_{t\to\infty}d(\f(sts'),\f(t))=0.
\]
Then we say that $(\Ga,\Lm)$ is \textit{bi-invariant at infinity relative to $X$}.

\medskip
The following lemma is taken from the proof of Theorem 3 in \cite{Oz1}; we reproduce its proof for the reader's convenience.

\begin{lemma}
Let $\Ga,\Lm$ and $X$ be as above, and assume that $(\Ga,\Lm)$ is bi-invariant at infinity relative to $X$. Then there exists a map $\zeta:\Lm\rightarrow \Pb(\Ga)$ such that, for all $g\in\Ga, s,s'\in \Lm$, 
\[
\lim_{y\to\infty}(\Vert g.\zeta(y)-\zeta(\alpha_g(y))\Vert_1+\Vert \zeta(sys')-\zeta(y)\Vert_1)=0.
\]
\end{lemma}
\begin{proof}
Since the action $\Ga\curvearrowright X$ is amenable, choose a sequence of continuous maps $\mu_n:X\rightarrow \Pb(\Ga)$ for $n\geq 1$ such that
\[
\lim_{n\to\infty}\Big(\sup_{x\in X}\|g.\mu_n^x-\mu_n^{g.x}\|_1\Big)=0\quad (g\in \Ga).
\]
Then define the map $\zeta_n:\Lm\setminus\{e\}\rightarrow \Pb(\Ga)$ by $\zeta_n(s)=\mu_n^{\f(s)}, s\in \Lm\setminus\{e\}$.

Taking a subsequence if necessary, we assume furthermore that
\[
\limsup_{t\to\infty}(\|\zeta_n(\alpha_g(t))-g.\zeta_n(t)\|_1+\|\zeta_n(sts')-\zeta_n(t)\|_1)<\frac{1}{n}
\]
for all $n,g\in\Ga, s,s'\in\Lm$. Indeed, for $g\in\Ga,s,s'\in\Lm$ one has:
\begin{multline*}
\|\zeta_n(\alpha_g(t))-g.\zeta_n(t)\|_1+\|\zeta_n(sts')-\zeta_n(t)\|_1=\\
\|\mu_n^{g.\f(t)}-g.\mu_n^{\f(t)}\|_1+\|\mu_n^{\f(sts')}-\mu_n^{\f(t)}\|_1\to_{t\to\infty}0.
\end{multline*}
For $n\geq 1,g\in\Ga,s,s'\in\Lm$, let us define the finite set
\[
D_n(g;s,s')\coloneqq \{t\in\Lm\setminus\{e\}\colon 
\|\zeta_n(\alpha_g(t))-g.\zeta_n(t)\|_1+\|\zeta_n(sts')-\zeta_n(t)\|_1\geq 1/n\}.
\]
Choose sequences of finite, symmetric sets $\{e\}=E_0\subset E_1\subset E_2\subset \cdots\subset \Ga$, $\bigcup_nE_n=\Ga$, and
$\{e\}=F_0\subset F_1\subset F_2\subset\cdots\subset \Lm$, $\bigcup_n F_n=\Lm$. Define inductively finite sets $\{e\}=\Omega_0\subset \Omega_1\subset \Omega_2\subset \cdots\subset \Lm$ by 
\[
\Omega_n=\bigcup_{g\in E_n,s,s'\in F_n,t\in\Omega_{n-1}}(D_n(g;s,s')\cup\{\alpha_g(t),sts'\})
\]
for each $n\geq 1$. For $y\in\Lm,$ set $l(y)=\min\{n\colon y\in \Omega_n\}$. Then finally, define $\zeta:\Lm\rightarrow \Pb(\Ga)$ by
\[
\zeta(y)=\frac{1}{l(y)}\sum_{n=0}^{l(y)-1}\zeta_n(y)
\]
 if $y\not=e$, and by $\zeta(y)=\delta_e$ if $y=e$.

Fix $g\in\Ga, s,s'\in\Lm$. Let $k$ be large enough so that $g\in E_k, s,s'\in F_k$. If $y\in\Lm$ is such that $l(y)>k$, then $y\in \Omega_{l(y)}$, but $y\notin \Omega_n$ for $n\leq k$. One has $\alpha_g(y),sys'\in \Omega_{l(y)+1}$ hence $|l(\alpha_g(y))-l(y)|,|l(sys')-l(y)|\leq 1$ and thus 
\[
\|\zeta_n(\alpha_g(y))-g.\zeta_n(y)\|_1+\|\zeta_n(sys')-\zeta_n(y)\|_1<1/n
\]
for $k\leq n<l(y)$. This proves that
\[
\lim_{l(y)\to\infty}(\|g.\zeta(y)-\zeta(\alpha_g(y))\|_1+\|\zeta(sys')-\zeta(y)\|_1)=0
\]
for all $g\in\Ga, s,s'\in\Lm$.
\end{proof}

The following lemma was suggested by A. Valette.

\begin{lemma}
For $n\geq 2$, let $\Pp^{n-1}(\R)$ be the real projective space of dimension $n-1$ equipped with the metric
\[
d([v],[w])=\arccos\Big(\frac{|\la v|w\ra|}{\Vert v\Vert \Vert w\Vert}\Big)\quad (v,w\in \R^n\setminus \{0\}).
\]
(For $v\in \R^n\setminus \{0\}$, $[v]$ denotes the image of $v$ by the canonical projection from $\R^n\setminus \{0\}$ onto $\Pp^{n-1}(\R)$.)
Then the map $\f:\Z^n\setminus \{0\}\rightarrow \Pp^{n-1}(\R)$ defined by $\f(v)=[v]$ is $\Ga_n$-equivariant, the action of $\Ga_n$ on $\Pp^{n-1}(\R)$ is amenable and one has for every $a\in \Z^n$
\[
\lim_{\Vert v\Vert\to\infty}d(\f(v+a),\f(v))=0.
\]
In particular, the pair $(\Ga_n,\Z^n)$ is bi-invariant at infinity relative to $\Pp^{n-1}(\R)$.
\end{lemma}
\begin{proof}
It is well-known that the linear action of $\SL_n(\R)$ on $\R^n$ induces an action on $\Pp^{n-1}(\R)$, and it was observed in the proof of Theorem 1.1 in \cite{JV} that $\rho_{n-1}(\SL_2(\Z))$ is a subgroup of $\SL_n(\Z)$. Thus, the restriction $\f$ to $\Z^n\setminus\{0\}$ of the canonical map $v\mapsto [v]$ to $\Pp^{n-1}(\R)$ is $\Ga_n$-equivariant. By \cite[Proposition 3.5]{JV}, the action of $\Ga_n$ on $\Pp^{n-1}(\R)$ is amenable, and one obviously has
\[
\lim_{\Vert v\Vert\to\infty} d([v+a],[v])=0
\]
for every $a\in \Z^n$ by definition of $d$. This proves that the pair $(\Ga_n,\Z^n)$ is bi-invariant at infinity. 
\end{proof}

The previous lemma implies that Theorem 1.1 is a special case of Theorem 1.2. Thus we are left to prove the latter.

\medskip\noindent
\textit{Proof of Theorem 1.2.}

We come back to the general situation where $\Ga$ is bi-exact, that it acts on the amenable group $\Lm$ and that there is a map $\zeta:\Lm\rightarrow \Pb(\Ga)$ such that, for all $g\in \Ga,s,s'\in\Lm$,
\[
\lim_{y\to\infty}(\Vert g.\zeta(y)-\zeta(\alpha_g(y))\Vert_1+\Vert \zeta(sys')-\zeta(y)\Vert_1)=0.
\]
We observe that, as exactness of groups is preserved under semi-direct products, the group $G\coloneqq\Lm\rtimes \Ga$ is exact.

In order to prove Theorem 1.2, by Theorem 2.2, it suffices to show the existence of a ucp map $\Phi:C^*_\lm(G)\otimes C^*_\rho(G)\rightarrow \B(\ell^2(G))$ such that $\Phi(a\otimes b)-ab\in \K(\ell^2(G))$ for all $a\in C^*_\lm(G),b\in C^*_\rho(G)$. We adapt Section 4 of \cite{Oz3} to our situation.
\medskip

Set $A=C^*_\lm(\Lm)$ and denote also by $\alpha$ the action of $\Ga$ by $*$-automorphisms on $A$ so that $\alpha_g(\lm(s))\coloneqq \lm(\alpha_g(s)), g\in \Ga,s\in\Lm$. Let $\pi_\alpha:A\rightarrow \B(\ell^2(\Lm)\otimes\ell^2(\Ga))$ be the (faithful) $*$-representation defined by
\[
\pi_\alpha(a)(\eta\otimes\delta_g)=\alpha_{g^{-1}}(a)\eta\otimes\delta_g \quad (a\in A,\eta\in\ell^2(\Lm),g\in\Ga).
\]
Then the pair $(\pi_\alpha,1\otimes\lm_\Ga)$ is a covariant representation of the $\Ga$-$C^*$-dynamical system $(A,\Ga,\alpha)$, i.e.
\[
1\otimes\lm_\Ga(g)\pi_\alpha(a)1\otimes\lm_\Ga(g)=\pi_\alpha(\alpha_g(a)) \quad (g\in\Ga,a\in A).
\]
Thus the reduced crossed product $A\rtimes_r\Ga\eqqcolon B$ is the $C^*$-algebra generated by $\pi_\alpha(A)\cup 1\otimes\lm_\Ga(\Ga)$,
and it is naturally $*$-isomorphic to $C^*_\lm(\Lm\rtimes \Ga)$. Hence it is locally reflexive  by \cite[Theorem 9.3.1]{BO} since $G=\Lm\rtimes\Ga$ is exact.

The unit vector $\delta_e\otimes\delta_e$ is tracial, cyclic and separating for $B$, and if the antilinear involution $J$ on $\ell^2(\Lm)\otimes\ell^2(\Ga)$ is defined by $J(x(\delta_e\otimes\delta_e))=x^*(\delta_e\otimes\delta_e)$ for every $x\in B$, then $C\coloneqq JBJ$ is included into the commutant of $B$.

Let also $u:\Ga\rightarrow U(\ell^2(\Lm))$ be the unitary representation associated to the action $\Ga\stackrel{\alpha}{\curvearrowright}\Lm$:
\[
(u(g)\eta)(s)=\eta(\alpha_{g^{-1}}(s))\quad (g\in \Ga, \eta\in\ell^2(\Lm),s\in\Lm).
\]
Then the following formula will be used below: for all $g\in\Ga,t\in\Lm$, one has
\[
u(g)\delta_t=\delta_{\alpha_g(t)}
\]
where $(\delta_t)_{t\in\Lm}$ denotes the natural orthonormal basis of $\ell^2(\Lm)$.

Set finally $\xi=\zeta^{1/2}$ and define $V:\ell^2(\Lm)\rightarrow \ell^2(\Lm)\otimes \ell^2(\Ga)=\ell^2(G)$ by $V\delta_y=\delta_y\otimes \xi(y)$ for every $y\in \Lm$. It is easy to see that $V$ is an isometry.

Thus, the map $\theta: \B(\ell^2(\Lm)\otimes\ell^2(\Ga))\rightarrow \B(\ell^2(\Lm))$ defined by
\[
\theta(T)=V^*TV\quad (T\in \B(\ell^2(\Lm)\otimes\ell^2(\Ga)))
\]
is ucp.

\begin{lemma}
With hypotheses and notation above, the following equalities hold:
\begin{enumerate}
\item [(a)] $J1\otimes\lm(g)J=u(g)\otimes\lm(g)$ for every $g\in \Ga$.
\item [(b)] $J\pi_\alpha(a)J=J_AaJ_A\otimes 1$ for every $a\in A$, where $J_A$ is the canonical antilinear involution on $\ell^2(\Lm)$ characterise by $J_A(a\delta_e)=a^*\delta_e$ for $a\in A$.
\item [(c)] $\theta((u\otimes \lm)(g))\delta_y=\la g.\xi(y)|\xi(\alpha_g(y))\ra\delta_{\alpha_g(y)}\quad (g\in \Ga,y\in\Lm)$;
\item [(d)] $\theta(\lm(s)\otimes 1)\delta_y=\la \xi(y)|\xi(sy)\ra\delta_{sy}\quad (s,y\in\Lm)$;
\item [(e)] $\theta(\rho(s)\otimes 1)\delta_y=\la \xi(y)|\xi(ys^{-1}\ra\delta_{ys^{-1}}\quad (s,y\in\Lm)$.
\end{enumerate}
Moreover, 
\[
\theta((u\otimes\lm)(g))-u(g),\ \theta(\lm(s)\otimes 1)-\lm(s),\ \theta(\rho(s)\otimes 1)-\rho(s)\in \K(\ell^2(\Lm))
\]
for all $g\in \Ga,s\in\Lm$.
\end{lemma}
\begin{proof} (a) We are going to prove the following relation first:  
\[
\delta_t\otimes \delta_h=\pi_\alpha(\lm(\alpha_h(t)))1\otimes\lm(h)(\delta_e\otimes \delta_e)
\]
for all $t\in \Lm,h\in \Ga$. Indeed,
\[
\pi_\alpha(\lm(\alpha_h(t)))1\otimes \lm(h)(\delta_e\otimes\delta_e)=
\pi_\alpha(\lm(\alpha_h(t)))(\delta_e\otimes \delta_h)=\alpha_{h^{-1}}(\lm(\alpha_h(t)))\delta_e\otimes \delta_h=\delta_t\otimes\delta_h.
\]
Thus
\begin{equation}\label{pi_alpha}
  J\delta_t\otimes\delta_h=1\otimes\lm(h^{-1})\pi_\alpha(\lm(\alpha_h(t^{-1})))(\delta_e\otimes\delta_e) 
\end{equation}
and 
\begin{multline*}
J1\otimes\lm(g)J\delta_t\otimes\delta_h=J1\otimes\lm(gh^{-1})\pi_\alpha(\lm(\alpha_h(t^{-1})))(\delta_e\otimes\delta_e)=\\
\pi_\alpha(\lm(\alpha_h(t)))(\delta_e\otimes\delta_{hg^{-1}})=\delta_{\alpha_g(t)}\otimes\delta_{hg^{-1}}=u(g)\delta_t\otimes\rho(g)\delta_h.
\end{multline*}
(b) We still use equality (\ref{pi_alpha}):
\begin{multline*}
J\pi_\alpha(\lm(s))J\delta_t\otimes\delta_h=J\pi_\alpha(\lm(s))1\otimes \lm(h^{-1})\pi_\alpha(\lm(\alpha_h(t^{-1})))(\delta_e\otimes\delta_e)=\\
\pi_\alpha(\lm(\alpha_h(t)))1\otimes \lm(h)\pi_\alpha(\lm(s^{-1})(\delta_e\otimes\delta_e)=
\pi_\alpha(\lm(\alpha_h(t)))(\delta_{s^{-1}}\otimes \delta_h)=\\
\lm(t)\delta_{s^{-1}}\otimes\delta_h=\delta_{ts^{-1}}\otimes\delta_h=\rho(s)\otimes 1\delta_t\otimes\delta_h.
\end{multline*}
(c) For $g\in \Ga, y,t\in \Lm$, one has
\begin{align*}
\la \theta((u\otimes \lm)(g))\delta_y|\delta_t\ra
&=
\la V^*(u\otimes\lm)(g)V\delta_y|\delta_t\ra=
\la \delta_{\alpha_g(y)}\otimes g.\xi(y)|\delta_t\otimes\xi(t)\ra\\
&=
\la \delta_{\alpha_g(y)}|\delta_t\ra\la g.\xi(y)|\xi(t)\ra
= \la\la g.\xi(y)|\xi(\alpha_g(y))\ra\delta_{\alpha_g(y)}|\delta_t\ra.
\end{align*}
Verifications of equalities (d) and (e) are similar and left to the reader. \\
Finally, for $g\in \Ga$ and $y\in\Lm$, one gets from (c)
\[
(\theta((u\otimes\lm)(g))-u(g))\delta_y=u(g)\Big(\la g.\xi(y)|\xi(\alpha_g(y))\ra-1\Big)\delta_y
\]
and the function $y\mapsto \la 
g.\xi(y)|\xi(\alpha_g(y))\ra-1$ belongs to $c_0(\Lm)$ because
\begin{align*}
|\la g.\xi(y)|\xi(\alpha_g(y))\ra-1|
&=
|\la g.\xi(y)|\xi(\alpha_g(y))-g.\xi(y)\ra|\\
&\leq
\Vert g.\xi(y)\Vert_2\Vert \xi(\alpha_g(y))-g.\xi(y)\Vert_2\\
&\leq
\Vert \zeta(\alpha_g(y))-g.\zeta(y)\Vert_1^{1/2}\to 0
\end{align*}
as $y\to\infty$. Thus $\theta((u\otimes\lm)(g))-u(g)\in \K(\ell^2(\Lm))$.
The last two claims are proved similarly.
\end{proof} 

 Following Section 4 of \cite{Oz3}, let $K_1$ be the hereditary $C^*$-subalgebra of $\B(\ell^2(\Lm)\otimes \ell^2(\Ga))$ generated by $\K(\ell^2(\Lm))\otimes \ell^\infty(\Ga)$. Then $B$ and $C$ are in the multiplier of $K_1$ and hence $K\coloneqq K_1\cap C^*(B,C)$ is an ideal in $C^*(B,C)$. 

As the bicommutant $A''$ of $A=C^*_\lm(\Lm)$ is injective, Lemma 2.5 shows that \cite[Propositions 4.3 and 4.4]{Oz3} apply so that the homomorphism $\mu:B\odot C\rightarrow C^*(B,C)/K$ defined by
\[
\mu\Big(\sum_{k=1}^m a_k\otimes x_k\Big)=\sum_{k=1}^m a_kx_k+K
\]
is continuous with respect to the $\min$-norm. 

Next, set $I_1=\B(\ell^2(\Lm))\otimes \K(\ell^2(\Ga))$ and 
$I\coloneqq I_1\cap C^*(B,C)$.
Then, as in the proof of \cite[Corollary 4.5]{Oz3}, the homomorphism $\Phi:B\odot C\rightarrow C^*(B,C)/(I\cap K)$ defined by
\[
\Phi\Big(\sum_{k=1}^m a_k\otimes x_k\Big)=\sum_{k=1}^m a_kx_k+(I\cap K)
\]
is continuous with respect to the $\min$-norm and has a ucp lifting since $B$ is locally reflexive (see e.g. the proof of \cite[Lemma 5]{Oz2}). As $I_1\cap K_1=\K(\ell^2(\Lm)\otimes \ell^2(\Ga))$, the group $G=\Lm\rtimes \Ga$ is bi-exact by Theorem 2.2, and $L(G)$ is a solid von Neumann algebra. \\
This ends the proof of Theorem 1.2.

\bibliographystyle{plain}
\bibliography{refSolid}

\end{document}